\documentclass[10pt]{article}
\usepackage{amsfonts,amssymb,amsmath}
\allowdisplaybreaks[3]

\begin{document}

\title{Symmetries of Coincidence Site Lattices of Cubic Lattices}

\author{P. Zeiner\\
Institute for Theoretical Physics \& CMS, TU Wien,\\
Wiedner Hauptsra{\ss}e 8--10, 1040 Vienna, Austria}

\maketitle

\begin{abstract}
We consider the symmetries of coincidence site lattices of $3$--dimensional
cubic lattices. This includes the discussion of the symmetry groups and
the Bravais classes of the CSLs. We derive various criteria and necessary
conditions for symmetry operations of CSLs. They are used to obtain
a complete list of the symmetry
groups and the Bravais classes of those CSLs that are generated by a rotation
through the angle $\pi$.

\end{abstract}

\newcommand{\ncd}{\newcommand}
\def\N{\mathbb N}
\def\Q{\mathbb Q}
\def\R{\mathbb R}
\def\Z{\mathbb Z}
\def\BEQ#1\EEQ{\begin{equation}#1\end{equation}}
\def\BAL#1\EAL{\begin{align}#1\end{align}}
\ncd{\gvct}[1]{{\boldsymbol{#1}}}
\ncd{\bma}{\begin{pmatrix}}
\ncd{\ema}{\end{pmatrix}}
\ncd{\Eq}[1]{Eq.~(\ref{#1})}
\ncd{\privateremark}[1]{}
\newtheorem{theo}{Theorem}[section]
\newtheorem{dfn}{Definition}[section]
\newtheorem{lemma}[theo]{Lemma}

\section{Introduction}

Coincidence site lattices (CSL) are an important tool to characterize and
analyze the structure of grain boundaries in crystals~\cite{boll70,boll82}.
At grain boundaries, two lattices with different orientation meet and it is
thus natural to consider the intersection of these two lattices. Although
grain boundaries are two--dimensional objects it proves useful to investigate
the three--dimensional intersection of these two
lattices~\cite{boll70,boll82}. Usually those
grain boundaries are preferred for which there is a high coincidence of
lattice sites.

CSLs for three-dimensional cubic lattices have been investigated by various
authors e.g.~\cite{grim74,gribo74,ble81,grim84,baa97} and they are well
understood. In particular one knows the coincidence rotations and has a handy 
parameterization for them, one knows the coincidence index $\Sigma$,
the number of different CSLs for
a given $\Sigma$, the generating functions, etc. Grimmer also discusses
the equivalence classes~\cite{grim74} and some symmetry aspects of the
CSLs~\cite{grim76,grim73}. The latter is based on the computation of CSLs up
to $\Sigma=199$, but up to now systematic approaches to the symmetries of CSLs
are rare and often unknown.\cite{andr78}

In this paper we discuss the symmetries and the Bravais classes of CSLs.
We start with recalling the basic notions and properties of CSLs. We review
the equivalence classes of CSLs by using an approach that stresses the symmetry
properties of the CSLs and that will prove useful in the determination
of the symmetry groups of the CSLs. We then state some propositions
on symmetry elements of CSLs in general and specialize them for
three--dimensional cubic lattices. They enable us to determine the symmetry
groups and Bravais classes for all CSLs that are generated by a rotation
that is equivalent to a rotation through the angle $\pi$. Finally,
we list all possible
symmetry groups and Bravais classes for all three types of cubic lattices.

Let us fix some of the notations first. If $m$ and $n$ are integers, then
$m|n$ means that $m$ divides $n$. We shall use the following convention
for vectors: $3$--dimensional vectors will be characterized by an arrow,
e.g. $\vec r$, whereas boldface letters denote $4$--dimensional vectors
and quaternions, e.g. $\gvct{q}$. The corresponding
inner (scalar) products will be written
as $\vec q\cdot\vec r$ and $\langle \gvct{q} | \gvct{r} \rangle$,
respectively. Furthermore $:=$ means ``is defined by''.

Let $\gvct{L}\subseteq\R^n$ be an $n$-dimensional lattice and $R$ a rotation.
Then $\gvct{L}(R)=\gvct{L}\cap R\gvct{L}$ is called a  coincidence site
lattice~(CSL) if it is a sublattice\footnote{We call $\gvct{L}'$ a sublattice
of $\gvct{L}$ if $\gvct{L}'\subseteq \gvct{L}$ i.e. if $\gvct{L}'$ is a subset
of $\gvct{L}$.} of finite index of $\gvct{L}$,
the corresponding rotation is called a coincidence
rotation~\cite{baa97,baaap94}. The
coincidence index $\Sigma(R)$ is defined as the index of $\gvct{L}(R)$ in
$\gvct{L}$. By index we mean the group theoretical index of $\gvct{L}(R)$ in
$\gvct{L}$, where we view $\gvct{L}(R)$ and $\gvct{L}$ as additive groups.
Physically, the index $\Sigma(R)$ is the ratio of the volume of the
(primitive) unit cells
of the lattices $\gvct{L}(R)$ and $\gvct{L}$. 

In the following, we specify $\gvct{L}$ to be a cubic lattice;
in particular, we assume $\gvct{L}=\gvct{L}_p=\Z^3$, i.e. a primitive cubic
lattice. This is an important case since the results of the primitive
cubic case can be easily extended to the face centered and body centered case.
We will do this in the end. Thus for the moment $\gvct{L}=\Z^3$. Then one
can show that a rotation is a coincidence rotation if and only if
it is a orthogonal matrix with rational entries~\cite{gribo74,grim84,baa97}.

Now any proper rotation in three--dimensional space can be parameterized by
quaternions (Cayley's parameterization)~\cite{koecheng,hurw,val,hardy}:
\BEQ
R(\gvct{r})=\frac{1}{|\gvct{r}|^2}
\begin{pmatrix}
\kappa^2+\lambda^2-\mu^2-\nu^2 & -2\kappa\nu+2\lambda\mu &
2\kappa\mu+2\lambda\nu\\
2\kappa\nu+2\lambda\mu & \kappa^2-\lambda^2+\mu^2-\nu^2 &
-2\kappa\lambda+2\mu\nu\\
-2\kappa\mu+2\lambda\nu & 2\kappa\lambda+2\mu\nu &
\kappa^2-\lambda^2-\mu^2+\nu^2 
\end{pmatrix},
\EEQ
where $\gvct{r}=(\kappa,\lambda,\mu,\nu)$ and
$|\gvct{r}|^2=\kappa^2+\lambda^2+\mu^2+\nu^2$. For the ease of the reader we
have listed some typical rotations and their corresponding
quaternions in Table~\ref{tabquatrot}.
Thus the rational orthogonal matrices can be parameterized by integral
quaternions, i.e. by quaternions with integral coefficients
$\kappa,\lambda,\mu,\nu$. Note that we will call a quaternion an integer
quaternion if it is an integral quaternion $(\kappa,\lambda,\mu,\nu)$
or the sum of an integral quaternion with the quaternion $1/2(1,1,1,1)$.
We call an integral quaternion $\gvct{r}=(\kappa,\lambda,\mu,\nu)$ primitive
if the greatest common divisor of $\kappa,\lambda,\mu,\nu$ equals $1$. If
not stated otherwise every (integral) quaternion will be assumed to be a
primitive quaternion.

Crystallographers may not be familiar with quaternions. Loosely speaking they
are four dimensional vectors that can be multiplied in a nice way. You can
view the quaternion $\gvct{r}=(\kappa,\lambda,\mu,\nu)$ as the
\mbox{$2\times 2$}--matrix
\BAL
\kappa\sigma_0+\lambda\sigma_1+\mu\sigma_2+\nu\sigma_3,
\EAL
where
\BAL
\sigma_0&=\bma 1&0\\ 0& 1\ema & \sigma_1&=\bma 0& i\\ i &0 \ema \\
\sigma_2&=\bma 0& 1\\ -1 &0 \ema & \sigma_3&=\bma i&0\\ 0& -i\ema,
\EAL
which are just the well known Pauli matrices (up to a factor $i$). The inner
product of two quaternions
$\gvct{q}=(\alpha,\beta,\gamma,\delta)$ and
$\gvct{r}=(\kappa,\lambda,\mu,\nu)$ is just the ordinary inner product of
$\R^4$
\BAL
\langle \gvct{q} | \gvct{r} \rangle
:=\alpha\kappa+\beta\lambda+\gamma\mu+\delta\nu.
\EAL
In addition, for any quaternion $\gvct{r}$ we define the conjugated
quaternion by $\bar{\gvct{r}}:=(\kappa,-\lambda,-\mu,-\nu)$. For more details
we have to refer to the literature~\cite{koecheng,hurw,val}.
Although we will
use quaternions extensively in the following no knowledge of quaternions
is necessary to understand most of the results. Those who are not interested
in the mathematical details may skip the proofs and simply keep in mind that
quaternions are a nice way to parameterize $3$--dimensional rotations.

As a matter of fact $R(\gvct{r})$ describes a proper rotation (i.e.
$\det R(\gvct{r})=1$ for all $\gvct{r}$) with rotation
axis $\vec{v}_0=(\lambda,\mu,\nu)^t$ and the rotation angle $\varphi$ given
by
\BAL
\cos\varphi=\frac{\kappa^2-\lambda^2-\mu^2-\nu^2}{
\kappa^2+\lambda^2+\mu^2+\nu^2}.
\EAL
In particular $\varphi=\pi$ for $\kappa=0$.

\begin{table}[htb]
\begin{tabular}{|c|c|}
\hline
$\gvct{q}$ & $R(\gvct{q})$\\
\hline
$(1,0,0,0)$ & $1$\\
$(0,1,0,0)$ & $2\ \ x,0,0$\\
$(0,0,1,0)$ & $2\ \ 0,y,0$\\
$(1,1,0,0)$ & $4^+\ \ x,0,0$\\
$(0,1,1,0)$ & $2\ \ x,x,0$\\
$(0,1,-1,0)$ & $2\ \ x,\bar x,0$\\
$(0,0,1,1)$ & $2\ \ 0,y,y$\\
$(0,1,1,1)$ & $2\ \ x,x,x$\\
$(1,1,1,1)$ & $3^+\ \ x,x,x$\\
$(3,1,1,1)$ & $6^+\ \ x,x,x$\\
$(m,n,0,0)$ & $\phi=\arccos\frac{m^2-n^2}{m^2+n^2}$,  $[100]$ \\
$(m,n,n,0)$ & $\phi=\arccos\frac{m^2-2n^2}{m^2+2n^2}$,  $[110]$ \\
$(m,n,n,n)$ & $\phi=\arccos\frac{m^2-3n^2}{m^2+3n^2}$,  $[111]$ \\
$(\kappa,\lambda,\mu,\nu)$ &
$\phi=\arccos\frac{\kappa^2-\lambda^2-\mu^2-\nu^2}{
\kappa^2+\lambda^2+\mu^2+\nu^2}$,  $[\lambda\mu\nu]$ \\
\hline
\end{tabular}
\caption{\label{tabquatrot} Quaternions and their symmetry operations:
The table lists the symmetry operations for several typical quaternions.
If the rotation is a crystallographic one, the standard crystallographic
notation is used. For proper coincidence rotations (i.e. non-crystallographic)
we state the rotation angle $\phi$ and the lattice direction of the
rotation axis.}
\end{table}

One can show that the coincidence index is given by
$\Sigma(R(\gvct{r}))=|\gvct{r}|^2/2^\ell$, where $\ell$ is the maximal power
such that $2^\ell$ divides $|\gvct{r}|^2$
(see e.g.~\cite{gribo74,grim84,baa97}).

Finally we mention a nice representation of the lattice vectors of a CSL.
To this end we define the vectors
\BAL\label{ri}
\vec{r}^{(0)}&=\bma r_1\\ r_2 \\ r_3 \ema, &
\vec{r}^{(1)}&=\bma r_0\\ r_3 \\ -r_2 \ema, \nonumber\\
\vec{r}^{(2)}&=\bma -r_3 \\ r_0\\ r_1 \ema, &
\vec{r}^{(3)}&=\bma r_2 \\ -r_1\\ r_0 \ema.
\EAL
for a primitive quaternion $\gvct{r}=(r_0,r_1,r_2,r_3)$. 
These vectors $\vec r^{(i)}$ are lattice vectors of $\gvct{L}(R(\gvct{r}))$.
Furthermore they are linearly dependent, in particular
$r_0\vec r^{(0)}-r_1\vec r^{(1)}-r_2\vec r^{(2)}-r_3\vec r^{(3)}=\vec 0$.
On the other hand, the vectors $\vec r^{(i)}$ span $\R^3$, hence we can obtain
a basis if we appropriately choose three of them. Thus any lattice vector
$\vec v\in\gvct{L}(R(\gvct{r}))$ is a rational linear combination of them.
However, such a basis has two disadvantages. First, rational coefficients are
not so handy as integer ones. Secondly there is no general rule on how to
choose the three basis vectors and any choice would break the symmetry of
the setting. These disadvantages can be avoided if we express the lattice
vectors of $\gvct{L}(R(\gvct{r}))$ as integer combinations of more than
three lattice vectors, which necessarily are linearly dependent. The following
lemma states how this can be achieved and we will often refer to it later on: 
\begin{lemma}\label{cslbasis}
The CSL $\gvct{L}(R(\gvct{r}))$ with $\gvct{r}=(r_0,r_1,r_2,r_3)$ is the
$\Z$--span of the following vectors:
\begin{itemize}
\item $\vec r^{(0)},\vec r^{(1)},\vec r^{(2)},\vec r^{(3)}$
if $|\gvct{r}|^2$ is odd,
\item $\vec r^{(0)},\vec r^{(1)},\vec r^{(2)},\vec r^{(3)},
1/2\,(\vec r^{(0)}+\vec r^{(1)}+\vec r^{(2)}+\vec r^{(3)})$ if $|\gvct{r}|^2$ even but not
divisible by $4$,
\item $\vec r^{(0)}, 1/2\,(\vec r^{(0)}+\vec r^{(1)}), 1/2\,(\vec r^{(0)}+\vec r^{(2)}),
1/2\,(\vec r^{(0)}+\vec r^{(3)})$ if $|\gvct{r}|^2$ is divisible by $4$.
\end{itemize}
\end{lemma}

\emph{Proof:} The CSL $\gvct{L}(R(\gvct{r}))$ with $\gvct{r}=(r_0,r_1,r_2,r_3)$
contains the vectors $\vec{r}^{(0)},\vec{r}^{(1)},\vec{r}^{(2)},
\vec{r}^{(3)},$, see e.g.~\cite{baa97}. 
Hence the lattices $\gvct{L}_i(\gvct{r})$ generated by the vectors
$\vec{r}^{(j)}, j\ne i$
are sublattices of $\gvct{L}$ with index $r_i|\gvct{r}|^2$ if $r_i\ne 0$.
Let $\gvct{L}'(\gvct{r})$ be the $\Z$--span of $\vec r^{(i)}, i=0,\ldots, 3$.
Then
$\gvct{L}'(\gvct{r})$ is a superlattice of $\gvct{L}_i(\gvct{r})$ and a sublattice
of $\gvct{L}(\gvct{r}))$. In order to determine the index of $\gvct{L}'(\gvct{r})$
in $\gvct{L}$ we first observe
$r_0\vec r^{(0)}-r_1\vec r^{(1)}-r_2\vec r^{(2)}-r_3\vec r^{(3)}=\vec 0$.
Thus we can express
any vector $\vec r^{(i)}$ as a linear combination of $\vec{r}^{(j)}, j\ne i$
if $r_i\ne 0$. Since the $r_i$ are relatively prime, $c\vec r^{(i)}, c\in \Z$
is an integer combination of $\vec{r}^{(j)}, j\ne i$
if and only if $r_i\big|c$.
Hence $\gvct{L}_i(\gvct{r})$ is a sublattice of $\gvct{L}'(\gvct{r})$ with index
$r_i$ if $r_i\ne 0$, and hence the index of $\gvct{L}'(\gvct{r})$ in
$\gvct{L}$ is $|\gvct{r}|^2=2^\ell\Sigma(R(\gvct{r}))$. Furthermore
$\gvct{L}'(\gvct{r})$ is a sublattice of $\gvct{L}(R(\gvct{r}))$ with index
$|\gvct{r}|^2/\Sigma(R(\gvct{r}))=2^\ell$.
Thus
$\gvct{L}'(\gvct{r})=\gvct{L}(R(\gvct{r}))$ if $|\gvct{r}|^2$ is odd. If
$|\gvct{r}|^2$ is even
then $1/2\,(\vec r^{(0)}+\vec r^{(1)}+\vec r^{(2)}+\vec r^{(3)})$
is a vector of
$\gvct{L}(R(\gvct{r}))$, and if $4\big | |\gvct{r}|^2$ then $\gvct{L}(R(\gvct{r}))$
contains all the vectors  $1/2\,(v_i\pm v_j)$. Since none of these vectors
is contained in $\gvct{L}'(\gvct{r})$, the claim follows.\hfill $\Box$

It follows from this lemma that every vector of the CSL can be written
in the form $m_0\vec r^{(0)}+m_1\vec r^{(1)}+m_2\vec r^{(2)}+m_3\vec r^{(3)}$, where the
coefficients $m_i$ are integers or half integers, with the following
constraints:
If $|\gvct{r}|^2$ is odd, all $m_i$ must be integers, if $|\gvct{r}|^2$ is even
then the sum $\sum_{i=0}^3m_i$ must be an integer. An additional constraint
applies if $|\gvct{r}|^2$ is even but not divisible by $4$. In this case
$\gvct{r}$ has exactly two odd and two even components. If the even components
are $r_{i_1}$ and $r_{i_2}$ and the odd components are $r_{i_3}$ and
$r_{i_4}$, then $m_{i_1}+m_{i_2}$ (as well as $m_{i_3}+m_{i_4}$) has to be an
integer, too. If we write $\gvct{m}=(m_0,m_1,m_2,m_3)$ then these constraints
can be summarized as follows: $\gvct{m}$ has an even number of integral
components and the product $\langle \gvct{r} | \gvct{m}\rangle$ is an integer.
An equivalent condition is that $2|\gvct{m}|^2$ and
$\langle \gvct{r} | \gvct{m}\rangle$ are integers.
Note that these coefficients $m_i$
are not unique, since the vectors $\vec r^{(i)}$ are linearly dependent.

\section{Equivalence classes of CSLs}

Different coincidence rotations $R$ may generate the same CSL or CSLs that
are just rotated versions of each other, so that an appropriate notion
of equivalence is desirable. Let $G$ be the point group of the
lattice~$\gvct{L}$. Then $\gvct{L}(RQ)=\gvct{L}(R)$ for all $Q\in G$. Moreover
$Q'RQ$ generates a rotated copy of $\gvct{L}(R)$, namely
$\gvct{L}(Q'RQ)=Q'\gvct{L}(R)$ for all $Q',Q\in G$.
These lattices are usually considered equivalent, since they are in a
crystallographic equivalent orientation with respect to the lattice $\gvct{L}$.
So we say that
two coincidence rotations $R,R'$ are {\em equivalent} if there exist rotations
$Q',Q\in G$ such that $R'=Q'RQ$. Some authors (e.g. Grimmer and
Bollmann~\cite{grim74,grim76,boll82}) extend this definition
and consider the lattices $\gvct{L}(R)$ and $\gvct{L}(R^{-1})$ as equivalent, too.
This is well justified from the physicist's point of view, since
$\gvct{L}(R^{-1})=R^{-1}\gvct{L}(R)$ is a rotated copy of $\gvct{L}(R)$.
Nevertheless $\gvct{L}(R^{-1})$ and $\gvct{L}(R)$ are in general not in a
crystallographically equivalent orientation with respect to $\gvct{L}$,
so we will use the more restrictive notion mentioned above throughout this
paper. Note that both definitions coincide if $R^2=1$, i.e. if $R$ is
a twofold rotation or a mirror reflection. 

It follows directly from the definition
that the set of all coincidence rotations equivalent to $R$ is the
double coset $GRG$, i.e. the set of all rotations $QRQ'$ with $Q,Q'\in G$.
The determination of all equivalence classes of
coincidence rotations is thus equivalent to the double coset decomposition
of the group of all coincidence rotations $OC(\gvct{L})$ with respect to the
subgroup $G$. Let $H=H(R)=G\cap RGR^{-1}$. Then $HRG=RG$. If
$G=\bigcup_{i} Q_i H$ is the coset decomposition of $G$ with respect to its
subgroup $H$, then we can express the double coset $GRG$ in terms of ordinary
cosets $GRG=\bigcup_{i} Q_i RG$. Thus $GRG$ consists of $|G|/|H|$ cosets,
and hence the number
of coincidence rotations equivalent to $R$ is given by $|GRG|=|G|^2/|H|$
(see \cite{hall} for some more details on double cosets).
Note that $H(R)$ and $H(R')$ are conjugated subgroups of $G$ if
$R$ and $R'$ are equivalent, in particular we have $H(R')=QH(R)Q^{-1}$ if
$R'=QRQ'$.

We now turn to the three--dimensional cubic case. For simplicity, we restrict
our considerations to proper rotations, i.e. $\det(R)=1$ or
$R\in SOC(\gvct{L})$, where $SOC(\gvct{L})$ denotes the group of all
coincidence rotations $R$ with $\det(R)=1$.
 But this is no real restriction, since any improper
rotation $R'$ is equivalent to $R$ if and only if $\bar1R'$ is equivalent to $R$,
where $\bar1$ denotes the inversion. The group of all proper rotations leaving
$\gvct{L}$ invariant is $O=432$ which is of order 24, the corresponding set of primitive
quaternions consists of the 48 quaternions $(\pm 1,0,0,0)$,
$(\pm 1,\pm 1,0,0)$, $(\pm 1,\pm 1,\pm 1,\pm 1)$ and permutations thereof.
If these quaternions are normalized to unity, they form a group, too,
namely the usual double cover\footnote{i.e. $O$ is a homomorphic image of the
corresponding group of quaternions and to each rotation there correspond
exactly two normalized quaternions, i.e. $\gvct{q}$ and $-\gvct{q}$ correspond
to the same rotation. This situation is well known from quantum mechanics,
where a rotation through $2\pi$ changes the sign of the spinor and only a
rotation about $4\pi$ leaves the spinor unchanged.} of $O$.

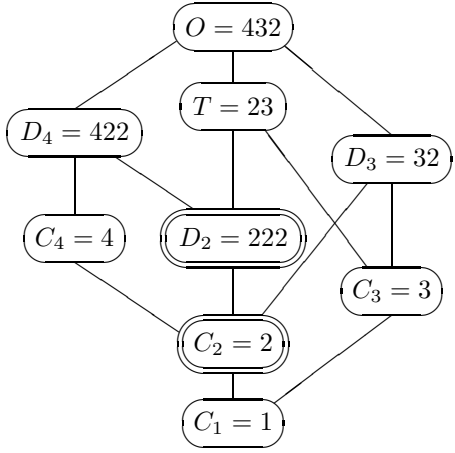
\begin{figure}

\unitlength10pt
\begin{picture}(20,19)
\put(11,17){\oval(4.5,1.8)}
\put(11,17){\makebox(0,0){$O=432$}}

\put(11,14.9){\line(0,1){1.2}}

\put(11,14){\oval(4.0,1.8)}
\put(11,14){\makebox(0,0){$T=23$}}

\put(5.0,13.9){\line(3,2){3.85}}

\put(5,13){\oval(5.1,1.8)}
\put(5,13){\makebox(0,0){$D_4=422$}}

\put(17,12.9){\line(-6,5){4.05}}

\put(17,12){\oval(4.5,1.8)}
\put(17,12){\makebox(0,0){$D_3=32$}}

\put(5,9.9){\line(0,1){2.2}}

\put(5,9){\oval(3.8,1.8)}
\put(5,9){\makebox(0,0){$C_4=4$}}

\put(5,8.1){\line(3,-2){3.95}}

\put(11,10.1){\line(0,1){3}}
\put(9.5,10.12){\line(-3,2){2.99}}

\put(11,9){\oval(5.5,2.2)}
\put(11,9){\oval(5.1,1.8)}
\put(11,9){\makebox(0,0){$D_2=222$}}

\put(16.1,7.9){\line(-3,4){3.9}}

\put(17,7.9){\line(0,1){3.2}}

\put(17,7){\oval(3.8,1.8)}
\put(17,7){\makebox(0,0){$C_3=3$}}

\put(11,6.1){\line(0,1){1.8}}

\put(12.1,6.12){\line(4,5){3.98}}

\put(11,5){\oval(4.2,2.2)}
\put(11,5){\oval(3.8,1.8)}
\put(11,5){\makebox(0,0){$C_2=2$}}

\put(11,2.9){\line(0,1){1}}

\put(17,6.1){\line(-4,-3){4.44}}

\put(11,2){\oval(3.8,1.8)}
\put(11,2){\makebox(0,0){$C_1=1$}}





\end{picture}

\caption{Subgroups of $432$. Single circles are used if all subgroups
of one type are conjugate under operations of $432$, double circles indicate
that the subgroups of these types are divided into two classes of conjugate
subgroups. For instance, there exist nine different subgroups of type $2$
corresponding to the nine rotations $2$ of $432$. Three of them are rotations
about the cubic fourfold axes and the other ones are rotations about the cubic
twofold axes. The fourfold axes are crystallographically equivalent 
directions, and so are the six twofold axes. Correspondingly we have two
classes of conjugate subgroups of type $2$, one class consisting of three
subgroups and the other one consisting of six subgroups.}
\label{subgroups}
\end{figure}

We have to determine the possible subgroups $H(R)\subset O$. Up to conjugacy,
$O$ has the following non trivial subgroups: the tetrahedral group $23$
of order $12$
generated by
\mbox{$3^+\ \ x,x,x$} and \mbox{$2\ \ x,0,0$},
the tetragonal group $422$ of order $8$
generated by \mbox{$4^+\ \ x,0,0$} and \mbox{$2\ \ 0,0,z$}, the tetragonal group $4$ of
order $4$ generated by \mbox{$4^+\ \ x,0,0$}, the trigonal group $32$
of order $6$ generated by \mbox{$3^+\ \ x,x,x$} and \mbox{$2\ \ x,\bar x,0$},
the trigonal group $3$ of order $3$
generated by \mbox{$3^+\ \ x,x,x$} and two orthorhombic subgroups $222$ of
order $4$. One has the generators \mbox{$2\ \ x,0,0$} and \mbox{$2\ \ 0,y,0$}
and the other one has \mbox{$2\ \ x,0,0$} and
\mbox{$2\ \ 0,y,y$}, respectively. Finally there exist two monoclinic groups
$2$ of order $2$
generated by \mbox{$2\ \ x,0,0$} and \mbox{$2\ \ x,x,0$}, respectively.
Not all of them can be realized in the form $H(R)=O\cap ROR^{-1}$, e.g.
the tetrahedral group is impossible. Fig.~\ref{subgroups} shows these
subgroups and the subgroup relations between them.

In order to determine the possible subgroups $H(R)$, it is convenient
to know the classes of conjugated elements of $O$. They are well known (see
e.g.~\cite{brad}) and are the following: the unit element,
the class of all rotations $3^{\pm}$ through $2\pi/3$,
the class 
of all rotations $4^{\pm}$ through $\pi/2$, and two classes of rotations through $\pi$,
namely $\{\mbox{$2\ \ x,0,0$};\ \mbox{$2\ \ 0,y,0$};\ \mbox{$2\ \ 0,0,z$}\}$
and
$\{\mbox{$2\ \ x,x,0$};\ \mbox{$2\ \ x,\bar x,0$};
\ \mbox{$2\ \ x,0,x$};\ \mbox{$2\ \ x,0,\bar x$};
\ \mbox{$2\ \ 0,y,y$};\ \mbox{$2\ \ 0,y,\bar y$}\}$.

Assume now that $H(R)$ contains the rotation $S=\mbox{$3^+\ \ x,x,x$}$. 
Due to the definition of $H(R)$ the threefold rotation $R^{-1}SR$ is contained
in $O$, too, i.e. it is a threefold rotation about the rotated axis
$R^{-1}\vec x$, $\vec x =(1,1,1)$, which is parallel to one of the
cubic threefold axis.
Since all threefold axis are crystallographically equivalent there exists
$Q\in O$ such that $R^{-1}\vec x=Q\vec x$, i.e. $R^{-1}SR=QSQ^{-1}$ or
$R^{-1}SR=QS^{-1}Q^{-1}$. In the latter case we make use of $S^{-1}=TST^{-1}$,
where $T=\mbox{$2\ \ x,\bar x,0$}$. We set $Q'=QT$ and get
$R^{-1}SR=Q'SQ'^{-1}$. Thus there always exists an appropriate $Q\in O$ such
that
$R^{-1}SR=QSQ^{-1}$, in fact this is just a consequence of the well known fact
that all threefold rotations are conjugate under operations of $O$.
Equivalently we may write 
$SRQ=RQS$, and hence $RQ$ commutes with
$S$. Now $RQ$ and $S$ can commute only if their rotation axes are parallel
and thus $RQ=R(m,n,n,n)$. Conversely $RQ=R(m,n,n,n)$ implies
$S\in H(RQ)=H(R)$. It immediately follows from geometric intuition
and it is straightforward to calculate that
$\mbox{$3^+\ \ \bar{x},x,x$}\in H(R)$ would imply
$RQ=R(1,1,1,1)=\mbox{$3^+\ \ x,x,x$}$ for an appropriate
$Q\in O$, and hence $H(R)=O$,
so that the tetrahedral group cannot be realized as $H(R)$. Similarly one
verifies that $\mbox{$2\ \ x,\bar x,0$}\in H(R)$ is possible if and only if
$m=0,\pm 1,\pm 3$ and $n=\pm 1$. This concludes the trigonal case.

One can proceed similarly if $S=\mbox{$4^+\ \ x,0,0$}\in H(R)$. Again $R$ must map
any rotation through $\pi/2$ onto a rotation of $\pi/2$, and since all these
rotations are conjugated elements in $O$, there exists a $Q\in O$ such that
$RQ$ commutes with $S$, and this statement holds if and only if $R$ is
equivalent to $R(m,n,0,0)$. One shows again that if $H(R)$ contains an
additional rotation about a twofold axis orthogonal to $(1,0,0)$, then
$H(R)=O$. Thus the only group $H(R)\ne O$ that contains $S$ is the
tetragonal group generated by $S=\mbox{$4^+\ \ x,0,0$}$.

It remains to check the orthorhombic and the monoclinic subgroups of $O$.
It turns out that only the monoclinic group $H(R)$ with 
generator \mbox{$2\ \ x,x,0$} can
be realized, which is the case if and only if $R$ is equivalent to
$R(m,n,n,0)$, where $0\ne |m|\ne|n|\ne 0$.
These observations can be summarized as follows:
\begin{theo}\label{equiv}
If $R$ is equivalent to the sixfold rotation $R(3,1,1,1)\sim R(0,1,1,1)$,
then the group $H(R)$ is conjugate to the trigonal group generated by
\mbox{$3^+\ \ x,x,x$} and \mbox{$2\ \ x,\bar x,0$} with $|H(R)|=6$. Thus there are $4\cdot 24$
(proper) rotations equivalent to $R(3,1,1,1)$.

If $R$ is equivalent to $R(m,n,n,n)$, $n\ne 0$, $m\ne\pm 3n$, then
$H(R)$ is conjugate to the trigonal group generated by \mbox{$3^+\ \ x,x,x$}
of order $|H(R)|=3$. There exist $8\cdot 24$ equivalent (proper) rotations.

If $R$ is equivalent to $R(m,n,0,0)$, $0\ne |m|\ne|n|\ne 0$, then
$H(R)$ is conjugate to the tetragonal group generated by
\mbox{$4^+\ \ x,0,0$} of order $|H(R)|=4$. There are $6\cdot 24$ equivalent (proper)
rotations.

If $R$ is equivalent to $R(m,n,n,0)$, $0\ne |m|\ne|n|\ne 0$, then $H(R)$
is conjugate to the monoclinic group generated by \mbox{$2\ \ x,x,0$}. Its
order $|H(R)|=2$ and thus there are $12\cdot 24$ equivalent (proper) rotations.

For all other $R\not\in O$, we have $|H(R)|=1$ and thus $24^2$ equivalent
(proper) rotations.
\end{theo}

Number theory provides explicit expressions for the number of inequivalent
rotations of the kind $(m,n,0,0)$, $(m,n,n,n)$, and $(m,n,n,0)$. The number
of representations of the binary forms $m^2+n^2$, $m^2+3n^2$, and $m^2+2n^2$
is well known and can be easily inferred from the prime decompositions
of $\Z[i]$, $\Z[e^{2\pi i/3}]$ and $\Z[i\sqrt{2}]$,
respectively~\cite{rib,hardy,dick}. The following theorem holds:
\begin{theo}
For given coincidence index $\Sigma$ there are $n_2$ inequivalent rotations
of the form $(m,n,0,0)$, $n_3$ rotations of the form $(m,n,n,n)$ and $n_4$
rotations of the form $(m,n,n,0)$. The total number $f_{ineq}(\Sigma)$
of inequivalent rotations is given by
\BAL
f_{ineq}(\Sigma)&= n_1 + n_2 + n_3 + n_4 + n_5,
\EAL
where
\BAL
n_1&=\left\{\begin{array}{rl}
1&\hspace{1cm}\mbox{if }\Sigma=3\\
0&\hspace{1cm}\mbox{otherwise}
\end{array}\right.\\
n_2&=\left\{\begin{array}{ll}
2^{m-1}&\hspace{0.5cm}\mbox{%
\begin{minipage}{5cm}\noindent
if $p=1\bmod 4$ for all prime factors $p$ of
$\Sigma$ and $m$ is the number of different prime factors of $\Sigma$.
\end{minipage}}\\
0&\hspace{0.5cm}\mbox{otherwise}
\end{array}\right.\\
n_3&=\left\{\begin{array}{ll}
2^{m-1}&\hspace{0.5cm}\mbox{%
\begin{minipage}{5cm}\noindent
if $p=1\bmod 6$ for all prime factors $p\ne 3$ of
$\Sigma>3$, the factor $p=3$ occurs at most once and $m$ is the number of
different prime factors $p=1\bmod 6$ of $\Sigma$.
\end{minipage}}\\
0&\hspace{0.5cm}\mbox{otherwise}
\end{array}\right.\\
n_4&=\left\{\begin{array}{ll}
2^{m-1}&\hspace{0.5cm}\mbox{%
\begin{minipage}{5cm}\noindent
if $p=1$ or $3\bmod 8$ for all prime factors $p$ of
$\Sigma$, where $m$ is the number of different prime factors of $\Sigma>3$.
\end{minipage}}\\
0&\hspace{0.5cm}\mbox{otherwise,}
\end{array}\right.
\EAL
and $n_5$ can be calculated from the total number of coincidence rotations
$24f(\Sigma)$,
\BAL
f(\Sigma)&= 4 n_1 + 6 n_2 + 8 n_3 + 12 n_4 + 24 n_5 \nonumber\\
&=\left\{\begin{array}{ll}0&\mbox{ if $\Sigma$ is even}\\
\Sigma\prod_{j}(1+\frac{1}{p_j})&\mbox{ if $\Sigma$ is odd},
\end{array}\right.
\EAL
where the product runs over all prime factors of $\Sigma$.
\end{theo}
For the last statement on the total number of coincidence rotations
see~\cite{baa97,grim76}. Note that $f(\Sigma)$ denotes the number of CSLs
of index $\Sigma$.

\section{Symmetries of CSLs}

\subsection{General remarks}

We turn now to the question which symmetries the CSLs have, i.e. we want to
find all rotations $Q$ such that $Q\gvct{L}(R)=\gvct{L}(R)$ holds. For the moment,
let $\gvct{L}$ be an arbitrary lattice. Then $Q\gvct{L}(R)=\gvct{L}(R)$
is certainly satisfied if $Q\in G\cap R G R^{-1}$, where $G$ denotes the
symmetry group of $\gvct{L}$. Furthermore $R$ is a symmetry operation if
$R^2\in G$, i.e. in particular if $R^2=1$. This follows immediately from
$R\gvct{L}(R)=R\gvct{L}\cap R^2\gvct{L}=R\gvct{L}\cap \gvct{L}=\gvct{L}(R)$
Thus, if $R$ is a rotation
through~$\pi$, then $R$ is a symmetry operation of $\gvct{L}(R)$. 
Thus a CSL $\gvct{L}(R)$ generated by a twofold rotation $R$ has at least
monoclinic symmetry.
More generally,
any rotation $Q\in RG$ such that $Q^2\in G$ is a symmetry operation
of $\gvct{L}(R)$. Let us define
\begin{dfn}\label{minsym}
The minimal symmetry group of $\gvct{L}(R)$ is the group generated by
$G\cap R G R^{-1}$ and all elements $Q\in RG$ such that
$Q^2\in G$.
\end{dfn}
It is clear that the minimal symmetry group is a subgroup of the symmetry group
of $G$. Naturally the question arises whether the symmetry group may contain
additional elements, and we will see below that the answer is affirmative.

Let $Q$ be a symmetry operation of $\gvct{L}(R)$. Then $Q\gvct{L}(R)=\gvct{L}(R)$
implies
\BEQ
\gvct{L}(R)=\gvct{L}\cap R\gvct{L}=\gvct{L}\cap R\gvct{L}\cap Q\gvct{L}\cap QR\gvct{L}
\subseteq \gvct{L}\cap Q\gvct{L} = \gvct{L}(Q).
\EEQ
Hence $|\gvct{L}(Q):\gvct{L}(R)|\in \N$ and therefore
$Q$ must be a coincidence rotation such that $\Sigma(Q)$ divides
$\Sigma(R)$. Similarly $\gvct{L}(R)\subseteq \gvct{L}(QR)$ and
$\Sigma(QR)\big|\Sigma(R)$. Conversely assume that
$\gvct{L}(R)\subseteq \gvct{L}(Q)$ and $\gvct{L}(R)\subseteq \gvct{L}(QR)$ hold.
Then $\gvct{L}(R)\subseteq Q\gvct{L}\cap QR\gvct{L}= Q\gvct{L}(R)$. Since $Q$
is orthogonal we must have $\gvct{L}(R)=Q\gvct{L}(R)$, and hence $Q$ is
a symmetry rotation of $\gvct{L}(R)$. We have thus proved
\begin{theo}\label{sym1}
$Q$ is a symmetry operation of $\gvct{L}(R)$ if and only if
$\gvct{L}(R)\subseteq \gvct{L}(Q)$ and $\gvct{L}(R)\subseteq \gvct{L}(QR)$ hold.
If $Q$ is a symmetry operation then $Q$ is a coincidence rotation
and $\gvct{L}(R)\subseteq \gvct{L}(Q^n R^i)$ for all $n\in\Z$ and $i=0,1$.
Thus $\Sigma(Q^n R^i)$ divides $\Sigma(R)$.
\end{theo}
The second part follows immediately from the first one since with $Q$ also
$Q^n$ is a symmetry operation.

\subsection{Symmetries of cubic lattices}

\subsubsection{Minimal symmetry groups}

The minimal symmetry groups for a three--dimensional cubic lattice follow
immediately from the preceding sections. Up to equivalence we have the
following cases:
\begin{itemize}
\item\label{minhex}
The minimal symmetry group for $R=R(3,1,1,1)$ is hexagonal and is generated
by $R(3,1,1,1)$ and $R(0,1,-1,0)$.
\item If $R$ is of the form $R(m,n,n,n)$, $0\ne m\ne\pm n,\pm 3n\ne 0$,
then the minimal symmetry group is trigonal and generated by
\mbox{$3^+\ \ x,x,x$} and $R(0,n+m,n-m,-2n)$.
\item\label{mintetr} If $R=R(m,n,0,0)$, $0\ne m \ne\pm n\ne 0$,
then the corresponding minimal symmetry group
is tetragonal. The generators are \mbox{$4^+\ \ x,0,0$} and
$R(0,0,m,n)$.
\item\label{minorth}
If $R=R(m,n,n,0)$, $0\ne m\ne\pm n\ne 0$
the minimal symmetry group is orthorhombic with
generators \mbox{$2\ \ x,x,0$} and $R(0,n,-n,m)$.
\item\label{minmon}
If $R=R(0,\ell,m,n)$ is not equivalent to one of the cases above,
then the minimal symmetry group is monoclinic and generated by $R$ itself.
\item If $R$ is not equivalent to a rotation through $\pi$, then the minimal
symmetry group is the trivial group consisting of the unit element only.
\end{itemize}

\subsubsection{Further symmetry operations}

Theorem~\ref{sym1} provides a criterion for symmetry operations.
In order to make use of it we need the following lemma.
\begin{lemma}\label{sym2}
Let $\gvct{q}$ and $\gvct{r}$ be primitive. Then
$\gvct{L}(R(\gvct{r}))\subseteq \gvct{L}(R(\gvct{q}))$ if there exists a
(half)integral quaternion $\gvct{m}$ such that $\gvct{r}=\gvct{q}\gvct{m}$.
\end{lemma}
\emph{Proof:}  
Lemma~\ref{cslbasis} provides us with a very convenient representation
of the lattice vectors of $\gvct{L}(R)$, which we will use very often in the
following. In \Eq{ri} we have defined the vectors $\vec r^{(i)}$ for a
quaternion $\gvct{r}$, analogously we define $\vec q^{(i)}$ for a
quaternion $\gvct{q}$.
 
Now $\gvct{r}=\gvct{q}\gvct{m}$ implies that
$\vec r=m_0 \vec q + m_1\vec{q}^{(1)} + m_2\vec{q}^{(2)} + m_3\vec{q}^{(3)}$,
and according to the comments after lemma~\ref{cslbasis} we have
$\vec r\in\gvct{L}(R(\gvct{q}))$ if both $2|\gvct{m}|^2$
and $\langle \gvct{q} | \gvct{m}\rangle$ are integers. But this
is certainly true since $r_0=\langle \gvct{q} | \bar{\gvct{m}}\rangle$ is an
integer.%
\privateremark{$2|\gvct{m}|^2$ must be an integer unless $|\gvct{r}|^2=1\pmod 4$
and $|\gvct{q}|^2=0 \pmod 4$. In the latter case all components are odd,
hence $r_0=\langle \gvct{q} | \bar{\gvct{m}}\rangle$ integer implies
that half integral components occur only pairwise.}
Similarly $\vec r^{(j)}\in \gvct{L}(R(\gvct{q}))$ follows from
$\gvct{r}\gvct{u_i}=\gvct{q}(\gvct{m}\gvct{u_i})$, where $\gvct{u_i}$
are the unit quaternions $(0,1,0,0)$, $(0,0,1,0)$, and $(0,0,0,1)$,
respectively. If $|\gvct{r}|^2$ is odd, this proves the lemma. If $|\gvct{r}|^2$
is even, then so is $|\gvct{q}|^2$ or $|\gvct{m}|^2$ (in the latter case
$\gvct{m}$ is integral).
Hence $\langle \gvct{q} | \gvct{n}\rangle$
and $2|\gvct{n}|^2$ are integers for
$\gvct{n}=1/2\,(\gvct{m}+\gvct{m}\gvct{u_1}+\gvct{m}\gvct{u_2}+\gvct{m}\gvct{u_3})$,
which has again only (half)integral components, and thus
$1/2\,(\vec r+\vec r^{(1)}+\vec r^{(2)}+\vec r^{(3)})\in\gvct{L}(R(\gvct{q}))$.
Similarly, one checks $1/2\,(\vec r^{(i)}+\vec r^{(j)})\in \gvct{L}(R(\gvct{q}))$
if $4$ divides $|\gvct{r}|^2$. Thus all generators of
$\gvct{L}(R(\gvct{r}))$ are in $\gvct{L}(R(\gvct{q}))$, which proves the
lemma.\hfill$\Box$

In order to determine the symmetry group of the CSL
$\gvct{L}(R), R=R(\gvct{r})$ it is
sufficient to determine all twofold rotations that leave $\gvct{L}(R)$
invariant. Thus we want to find all (primitive)
quaternions $\gvct{q}=(0,q_1,q_2,q_3)$
such that $Q=R(\gvct{q})$ is a symmetry operation of $\gvct{L}(R)$. In the
following, we shall always assume that $\gvct{q}$ is of the form
$\gvct{q}=(0,q_1,q_2,q_3)$.

First we state a simple result.
\begin{lemma}\label{sym3}
Let $\gvct{r}$ and $\gvct{q}=(0,q_1,q_2,q_3)$ be primitive quaternions. Then
$Q=R(\gvct{q})$ is a symmetry operation of $\gvct{L}(R(\gvct{r}))$ if there exists
a (half)integral quaternion $\gvct{m}$ such that $\gvct{r}=\gvct{q}\gvct{m}$ and
$\gvct{q}\gvct{m}=-\gvct{m}\gvct{q}$ or $\gvct{q}\gvct{m}=\gvct{m}\gvct{q}$.
\end{lemma}
\emph{Proof:}
According to lemma~\ref{sym2} we have $\gvct{L}(R)\subseteq\gvct{L}(Q)$.
Moreover $QR=R(\gvct{m})$, and a second application of lemma~\ref{sym3}
(maybe we have to shift a factor $2$ from $\gvct{q}$ to $\gvct{m}$)
gives $\gvct{L}(R)\subseteq\gvct{L}(QR)$, and hence by theorem~\ref{sym1}
$Q$ is a symmetry operation of $\gvct{L}(R)$.\hfill$\Box$

The importance of this lemma lies in the fact that this exhausts more or
less all cases for $\gvct{r}=(0,r_1,r_2,r_3)$. 
The condition $\gvct{q}\gvct{m}=-\gvct{m}\gvct{q}$
implies that the real part of $\gvct{q}$ and $\gvct{m}$ vanishes, i.e. that
$\gvct{q}$ and $\gvct{m}$ correspond to rotations through $\pi$,
moreover the two rotation axes are orthogonal to each other. Corresponding CSLs
have thus at least orthorhombic symmetry. The second case will be important in
deciding whether a CSL with minimal trigonal symmetry is hexagonal or not.

We want to find some necessary conditions for symmetry elements $Q$.
Since the vectors $\vec r^{(i)}$ are elements of $\gvct{L}(Q)$ by
theorem~\ref{sym1}, 
we infer from lemma~\ref{cslbasis} that they can be written as
$\vec r^{(i)}=n_0^{(i)}\vec q + n_1^{(i)}\vec q^{(1)} +
n_2^{(i)}\vec q^{(2)} + n_3^{(i)}\vec q^{(3)}$,
where the coefficients $n_j^{(i)}$ are integers or half-integers.
If we define $\hat R=(\vec r,\vec r^{(1)},\vec r^{(2)},\vec r^{(3)})$
we can reformulate these equations as $\hat R=\hat Q N$, where the entries
of the $4\times 4$ matrix $N$ are integers or half integers. But $\hat R$
can be obtained from $\pi(\gvct{r})$ by skipping the first row, where
$\pi$ is the matrix representation of the quaternions defined by
\BAL
\pi(1,0,0,0)&=\tau_0=
\bma
1 & 0 & 0 & 0 \\
0 & 1 & 0 & 0 \\
0 & 0 & 1 & 0 \\
0 & 0 & 0 & 1
\ema \\
\pi(0,1,0,0)&=\tau_1=
\bma
0 & -1 & 0 & 0 \\
1 & 0 & 0 & 0 \\
0 & 0 & 0 & -1 \\
0 & 0 & 1 & 0
\ema\\
\pi(0,0,1,0)&=\tau_2=
\bma
0 & 0 & -1 & 0 \\
0 & 0 & 0 & 1 \\
1 & 0 & 0 & 0 \\
0 & -1 & 0 & 0
\ema \\
\pi(0,0,0,1)&=\tau_3=
\bma
0 & 0 & 0 & -1 \\
0 & 0 & -1 & 0 \\
0 & 1 & 0 & 0 \\
1 & 0 & 0 & 0
\ema.
\EAL 
Thus $\pi(\gvct{r})=\pi(\gvct{q})N+X$ for an appropriately chosen matrix
\BAL
X=\bma
x_0 & x_1 & x_2 & x_3 \\
0 & 0 & 0 & 0 \\
0 & 0 & 0 & 0 \\
0 & 0 & 0 & 0
\ema.
\EAL
Hence there exists a (rational) quaternion $\gvct{m}$ such that
$N=\pi(\gvct{m})-\frac{1}{|\gvct{q}|^2}\pi(\bar{\gvct{q}})X$, explicitly
\BAL
\bma
n_0^{(0)} & n_0^{(1)} & n_0^{(2)} & n_0^{(3)}\\
n_1^{(0)} & n_1^{(1)} & n_1^{(2)} & n_1^{(3)}\\
n_2^{(0)} & n_2^{(1)} & n_2^{(2)} & n_2^{(3)}\\
n_3^{(0)} & n_3^{(1)} & n_3^{(2)} & n_3^{(3)}
\ema 
&= \bma
m_0 & -m_1 & -m_2 & -m_3 \\
m_1 & m_0 & -m_3 & m_2 \\
m_2 & m_3 & m_0 & -m_1 \\
m_3 & -m_2 & m_1 & m_0 
\ema 
+\frac{1}{|\gvct{q}|^2} \bma
0 & 0 & 0 & 0 \\
q_1 x_0 & q_1 x_1 & q_1 x_2 & q_1 x_3 \\
q_2 x_0 & q_2 x_1 & q_2 x_2 & q_2 x_3 \\
q_3 x_0 & q_3 x_1 & q_3 x_2 & q_3 x_3
\ema
\EAL
since $q_0=0$. Now $n_j^{(i)}$ are integer or half-integer and hence so are
the $m_i$. Thus $Q$ is a symmetry operation of $\gvct{L}(R)$ only if there exists
a (half) integral
quaternion $\gvct{m}$ such that $\gvct{r}=\gvct{q}\gvct{m}$.
Since $\gvct{r}$ is primitive, so is $2^\ell\gvct{m}$, where $\ell=0,1$
according to whether $\gvct{m}$ is integral or half integral. Thus we can state
\begin{lemma}\label{sym4}
The quaternion $\gvct{q}=(0,q_1,q_2,q_3)$ corresponds to a symmetry operation
of $\gvct{L}(R(\gvct{r}))$ only if $\gvct{m}:=1/|\gvct{q}|^2 \bar{\gvct{q}}\gvct{r}$
is a half integral quaternion.
\end{lemma}

Note in passing that $2x_i$ is a multiple of $|\gvct{q}|^2$ (since the
$q_i$ are relatively prime). Thus
$N'=N+\frac{1}{|\gvct{q}|^2}\pi(\bar{\gvct{q}})X=\pi(\gvct{m})$ is a (half)integer
matrix and satisfies $\hat R=\hat Q N'$, too (remember that $N$ is not
uniquely defined).

We knew from lemma~\ref{cslbasis} and theorem~\ref{sym1}
that $\vec r^{(i)}$ is an element of
$\gvct{L}(Q)$ and this has lead us to the representation
$\gvct{r}=\gvct{q}\gvct{m}$. But we know even more about $\vec r^{(i)}$.
If $Q$ is a symmetry operation of $\gvct{L}(R)$, then $Q\vec r^{(i)}$
and hence
\BAL
\vec r^{(i)}+Q\vec r^{(i)}=\frac{2\vec q\cdot\vec r^{(i)}}{|\vec q|^2} \vec q=
2\bar m_i \vec q
\EAL
are elements of $\gvct{L}(R)$, where $\bar m_i$ are the components of
the conjugate $\bar{\gvct{m}}$ of $\gvct{m}$.
%
%
Since $2^\ell\gvct{m}$ is primitive we infer $2^{1-\ell}\vec q\in \gvct{L}(R)$. 
We are now in a similar situation as before. However, we do not know
whether $2^{1-\ell}\vec q^{(i)}\in \gvct{L}(R)$ or not, 
so we can only conclude
that there exist a (half)integral quaternion $\gvct{n}$ and an integer
$x=\langle \gvct{r} | \gvct{n}\rangle$
such that $2^{1-\ell}\gvct{q}=\gvct{r}\gvct{n}-x\gvct{e}$, where $\gvct{e}$ is the
unity quaternion. Expressing $\gvct{q}$ in terms of $\gvct{r}$ and $\gvct{m}$
we get
\BAL\label{condsym1}
\frac{2^{1-\ell}}{|\gvct{m}|^2}\bar{\gvct{m}}
=\gvct{n}-\frac{x}{|\gvct{r}|^2}\bar{\gvct{r}}.
\EAL
This equation can hold only if $|\gvct{q}|^2$ divides
$2x=2\langle \gvct{r} | \gvct{n}\rangle$, hence
$y:=\langle \gvct{r} | \gvct{n}\rangle/|\gvct{q}|^2$ is integer or half integer.
Now
\BAL
|\gvct{m}|^2 |\gvct{n}|^2= 2^{2-2\ell}+y^2|\gvct{q}|^2
\EAL
since $\langle \gvct{r} | \gvct{m}\rangle=0$ due to $q_0=0$. But this implies
that the greatest common divisor (=$\gcd$) of $2|\gvct{m}|^2$
and $|\gvct{q}|^2$
is a power of $2$, and so is
and $\gcd(2|\gvct{m}|^2, 2y)$. Hence $\Sigma(\gvct{q})$ and $\Sigma(\gvct{m})$
are relatively prime. Thus we have proved:
\begin{theo}
A twofold rotation $Q=R(\gvct{q})$ can be a symmetry operation of
$\gvct{L}(R(\gvct{r}))$ only if there exists a (half) integral quaternion
$\gvct{m}$ such that $\gvct{r}=\gvct{q}\gvct{m}$ and $\Sigma(\gvct{q})$ and
$\Sigma(\gvct{m})$ are relatively prime. In particular, if $\Sigma(\gvct{r})$
is a prime power, then the symmetry group of $\gvct{L}(R(\gvct{r}))$ is just
the minimal symmetry group.
\end{theo}

Next we want to have a closer look on \Eq{condsym1}. Since it is too difficult
to discuss this equation in full generality we start with the case that
$R(\gvct{r})$ is a twofold rotation, i.e. $\gvct{r}=(0,r_1,r_2,r_3)$.
Let us further assume that $|\gvct{q}|^2$ and $|\gvct{r}|^2$ are both
odd, hence $\gvct{m}$ and $\gvct{n}$ are integral quaternions. Thus
\BAL
2\gvct{m}=-y\gvct{r} \bmod |\gvct{m}|^2
\EAL
and in particular $2m_0=0\bmod |\gvct{m}|^2$. Hence $m_0=0$ or
$2|m_0|\geq m_0^2$. The latter is satisfied only for the quaternions
$(\pm 2,0,0,0)\sim (\pm 1,0,0,0)$, $(\pm 1,\pm 1,0,0)$ and permutations
thereof. Since
$\Sigma(\gvct{m})=1$ for all of them, these are trivial solutions
(if they are solutions at all), hence it
remains to consider $m_0=0$. In this case $\gvct{m}$ corresponds to a rotation
through $\pi$, too, and the rotation axes $\vec r$, $\vec q$ and $\vec m$ must
be mutually orthogonal. In particular the equations $\gvct{r}=\gvct{q}\gvct{m}$
and $\gvct{q}=1/|\gvct{m}|^2\gvct{r}\bar{\gvct{m}}$ are equivalent to
$\vec r=\vec q\times \vec m$ and $\vec q=(1/\vec m^2) \vec m\times\vec r$.

Similarly we can handle the case $|\gvct{r}|^2$ odd and $|\gvct{q}|^2$ even.
Now $\gvct{m}=1/2\,\gvct{m'}$, $\gvct{m'}$ integral and $|\gvct{m'}|^2$ even.
$\gvct{n}$ is again integral. Thus \Eq{condsym1} reads
\BAL
\gvct{m'}=-2y\gvct{r} \bmod |\gvct{m'}|^2/2
\EAL
and we get the same trivial solutions (but now for $\gvct{m'}$) as before,
unless $2m_0=m_0'=0$. Thus again $\vec r$, $\vec q$ and $\vec m$ are
orthogonal.

The case $|\gvct{r}|^2$ even is slightly more difficult since $\gvct{n}$
may be half integral now. Again we assume $|\gvct{q}|^2$ odd first, hence
$\gvct{m}$ is integral and $|\gvct{m}|^2$ is even but
$4\!\!\not\big|\, |\gvct{m}|^2$.
We now replace
condition~(\ref{condsym1}) by the weaker condition
\BAL
2\gvct{m}=-y\gvct{r} \bmod |\gvct{m}|^2/2.
\EAL
The only possible solutions for $\gvct{m}$ with $m_0\ne 0$ such that $2$
but not $4$ divides $|\gvct{m}|^2$ are $(\pm 1,\pm 1,0,0)$ and permutations
thereof, hence non trivial solutions are again only possible for $m_0=0$,
and as before $\vec r$, $\vec q$ and $\vec m$ must be mutually orthogonal.

The last case to consider is $|\gvct{r}|^2$ and $|\gvct{q}|^2$ both even.
Then $|\gvct{m}|^2$ is an odd integer, but note that $\gvct{m}$ may be half
integral. Instead of (\ref{condsym1}) we consider again a weaker condition
\BAL
2\gvct{m}=-2y\gvct{r} \bmod |\gvct{m}|^2.
\EAL
Everything is the same as before, except that 
here $\gvct{m}=(\pm\frac{3}{2},\pm\frac{1}{2},\pm\frac{1}{2},\pm\frac{1}{2})$
are possible solutions. Of course it can only be a solution if
$3\big ||\gvct{r}|^2$. Assume that
$\gvct{m}=(\frac{3}{2},\frac{1}{2},\frac{1}{2},\frac{1}{2})$, which describes
a rotation through the angle $\pi/3$, is a solution. Due to
$\langle \gvct{r} | \gvct{m}\rangle=0$ we have
$\gvct{r}=(0,r_1,r_2,-r_1-r_2)$, which is equivalent to
$(r_1-r_2,r_1+r_2,r_1+r_2,r_1+r_2)$.

Note that we have proved so far only that
$\gvct{m}=(\frac{3}{2},\frac{1}{2},\frac{1}{2},\frac{1}{2})$ may be a solution
of \Eq{condsym1} if (and only if) $\gvct{r}=(0,r_1,r_2,-r_1-r_2)$ and
$3\big ||\gvct{r}|^2=2(r_1^2+r_2^2+r_1r_2)$.
We have still to prove that \Eq{condsym1} is indeed
satisfied. Inserting $\gvct{r}$ and $\gvct{m}$ in \Eq{condsym1} we obtain
\BAL
\gvct{n}&=\frac{1}{|\gvct{m}|^2}\bar{\gvct{m}}+\frac{x}{|\gvct{r}|^2}\bar{\gvct{r}}
=\frac{1}{6}(3,1,1,1)+
\frac{x}{2(r_1^2+r_2^2+r_1r_2)}(0,r_1, r_2, -r_1-r_2)
\EAL
and this is in fact a half integral quaternion if we choose
$x=2 r_1(r_1^2+r_2^2+r_1r_2)/3$ and take into account that
$r_1=r_2\ne 0\pmod 3$.

We can thus formulate the following lemma:
\begin{lemma}
Let $\gvct{r}=(0,r_1,r_2,r_3)$ and $\gvct{q}=(0,q_1,q_2,q_3)$ be primitive.
Then $Q=R(\gvct{q})$ is a symmetry operation of $\gvct{L}(R(\gvct{r}))$ if and
only if one of the following conditions hold:
\begin{itemize}
\item $Q$ is an element of the minimal symmetry group.
\item There exists a (half)integral quaternion $\gvct{m}=(0,m_1,m_2,m_3)$
such that $\gvct{r}=\gvct{q}\gvct{m}$, i.e. there exists a (half)integral vector
$\vec m$ such that $\vec r=\vec q\times \vec m$ and $\vec q\cdot\vec m=0$.
\item $\gvct{r}$ is equivalent to
$\gvct{r}'=\gvct{u}\gvct{r}\gvct{u}^{-1}=(0,r_1,r_2,-r_1-r_2)$ and
$\gvct{r}'=1/2\, \gvct{u}\gvct{q}\gvct{u}^{-1}(3,1,1,1)$, i.e.
$3\big | |\gvct{r}|^2$ and
$\gvct{q}=1/6\,\gvct{r}\gvct{u}^{-1}(3,-1,-1,-1)\gvct{u}$,
where $\Sigma(\gvct{u})=1$.
\end{itemize}
\end{lemma}
\emph{Proof:}
We have already proved that it is necessary that one of these conditions
holds. The first condition is sufficient by def.~\ref{minsym}.
From lemma~\ref{sym3} it follows immediately that the second condition
is sufficient, too. Thus it remains to show that
the third condition is sufficient. Let
$\gvct{r}=(0,r_1,r_2,-r_1-r_2)$ and $3\big | |\gvct{r}|^2$. Then
$r_1=r_2\pmod 3$ and hence
$\gvct{q}=1/6\,\gvct{r}(3,-1,-1,-1)=1/3(0,r_1-r_2,r_1+2r_2,-2r_1-r_2)$ has
integral components. Moreover $\gvct{r}$ is equivalent to
$1/2\,\gvct{r}(1,-1,-1,-1)=1/2\,\gvct{q}(0,1,1,1)=-1/2\,(0,1,1,1)\gvct{q}$,
and again lemma~\ref{sym3} concludes the proof.\hfill$\Box$

Note that $\gvct{r}=(0,r_1,r_2,-r_1-r_2)$ is equivalent to a quaternion
of the form $(m,n,n,n)$, and the last statement is equivalent to the fact
that $R(0,1,1,1)$ and hence $R(3,1,1,1)$ are symmetry elements of
the corresponding CSL if and only if $3$ divides $\Sigma(R(m,n,n,n))$, i.e.
$\gvct{L}(R(m,n,n,n))$ has hexagonal symmetry if and only if $3$ divides
$\Sigma(R(m,n,n,n))$.

We can go a step further and ask how many ways there are to write
$\gvct{r}$ as a product of two vectorial quaternions $\gvct{q}=(0,\vec q)$ and
$\gvct{m}=(0,\vec m)$
with $\Sigma(\gvct{q})$ and $\Sigma(\gvct{m})$ relatively prime.
i.e., we want to determine the number of twofold rotation axes orthogonal
to $\vec r$. Obviously to each decomposition $\gvct{r}=\gvct{q}\gvct{m}$ there
correspond two twofold rotation axes, namely $\vec q$ and $\vec m$, and to
each decomposition $\gvct{r}=\gvct{q}\gvct{m}$ there correspond the three
additional decompositions $\gvct{r}=(-\gvct{q})(-\gvct{m})=-\gvct{m}\gvct{q}
=\gvct{m}(-\gvct{q})$. If there are further decompositions, then there must
exist additional twofold rotation axes orthogonal to $\vec r$. But this is
impossible unless $\vec r$ is a fourfold or a sixfold rotation axis.
If $\vec r$ is a fourfold axis, then $\Sigma(\gvct{m})=\Sigma(\gvct{q})$, and
hence $\Sigma(\gvct{m})=\Sigma(\gvct{q})=\Sigma(\gvct{r})=1$, i.e. $\gvct{r}$
is equivalent to $(0,1,0,0)$. Similarly if $\vec r$ is a sixfold rotation
axis, then $\gvct{r}$ must be equivalent to $(0,1,1,1)$, since
$\Sigma(\gvct{q})=3$ and $\Sigma(\gvct{m})=1$ or vice versa. Thus, if
a decomposition $\gvct{r}=\gvct{q}\gvct{m}$ exists,
it is unique up to trivial operations except for the two cases mentioned
above.

With the knowledge we have developed so far we can immediately prove
\begin{theo}\label{symcsl}
Let $\gvct{r}$ be equivalent to a vectorial quaternion $(0,r_1,r_2,r_3)$.
Then the symmetry group of $\gvct{L}(R(\gvct{r}))$ is either the minimal symmetry
group or the minimal symmetry group is a subgroup of it of order $2$.
In particular,
we have (always $0\ne n\ne m\ne 0$)
\begin{enumerate}
\item If $\gvct{r}=(m,n,n,n)$ and $3$ does not
divide $\Sigma$, then the CSL has rhombohedral symmetry and its symmetry group
is just the minimal symmetry group, generated by \mbox{$3^+\ \ x,x,x$} and
$R(0,n+m,n-m,-2n)$. There
are precisely $n_3$ inequivalent CSLs for a fixed $\Sigma$.
\item\label{hex}
If $\gvct{r}=(m,n,n,n)$ and $3$ divides $\Sigma$,
then the CSL has hexagonal symmetry. If $\Sigma=3$ the symmetry group is the
minimal symmetry group generated by $R(3,1,1,1)$ and $R(0,1-1,0)$. If
$\Sigma>3$ the symmetry group is a proper supergroup of index $2$
of the minimal symmetry group and is generated by $R(3,1,1,1)$ and 
$R(0,n+m,n-m,-2n)$. There are again
$n_3$ inequivalent CSLs for a fixed $\Sigma$, except for $\Sigma=3$, where
we have only one.
\item\label{tetr} If $\gvct{r}=(m,n,0,0)$, the CSL has tetragonal symmetry,
its symmetry group is the minimal symmetry group generated by \mbox{$4^+\ \ x,0,0$} and
$R(0,0,m,n)$. There are $n_2$ inequivalent CSLs.
\item\label{orth1}
If $\gvct{r}=(m,n,n,0)$, the CSL has orthorhombic
symmetry. Its symmetry group is again just the minimal symmetry group
generated by \mbox{$2\ \ x,x,0$} and
$R(0,n,-n,m)$. There are $n_4$ inequivalent CSLs.
\item\label{orth2}
If $\gvct{r}=(0,\vec r)$ is not equivalent to one of the cases above,
then the CSL is orthorhombic if there exist two orthogonal integer vectors
$\vec q$ and $\vec m$, $\vec q^2$ and $\vec m^2$ relatively prime, such that
$\vec r=\vec q\times\vec m$ or $\vec r=1/2\ \vec q\times\vec m$. This
condition is equal to the existence of two integral quaternions $\gvct{q}$,
$\gvct{m}$ such that
$\gvct{r}=1/2^\ell \gvct{q}\gvct{m}=-1/2^\ell \gvct{m}\gvct{q}$.
If no such decomposition exists then the CSL is
monoclinic. The symmetry group of the former is generated by
$R(0,\vec r)$ and $R(0,\vec q)$, whereas the latter is generated by
$R(0,\vec r)$. If $\Sigma$ is a prime power, only the latter case is
possible.
\end{enumerate}
\end{theo}
These results are in coincidence with the observations of
W.~Grimmer~\cite{grim76}, who has calculated the CSLs and there symmetries
up to $\Sigma=199$.

This theorem covers all CSLs where $R$ is equivalent to a twofold operation.
For the general case, only some partial answers exist.
We want to discuss here only the question under what conditions
$(0,r_1,r_2,r_3)$ describes a symmetry operation of 
$\gvct{L}(R(\gvct{r}))$, $\gvct{r}=(r_0,r_1,r_2,r_3)$. We assume $r_0\neq 0$,
for otherwise the answer is trivial. Let $\gvct{q}=1/c (0,r_1,r_2,r_3)$,
where $c$ is the greatest common divisor of $r_1,r_2,r_3$. If $\gvct{q}$
describes a symmetry operation, then
$|\gvct{q}|^2$ must divide $2|\gvct{r}|^2$ and hence $2r_0^2$. From
lemma~\ref{sym4} we infer that
$1/|\gvct{q}|^2 \bar{\gvct{q}}\gvct{r}= 1/|\gvct{q}|^2 (c|\gvct{q}|^2,r_0\vec q)$
must be a half integral quaternion, hence $|\gvct{q}|^2$ must divide $2r_0$.
Conversely, assume that $|\gvct{q}|^2$ divides $2r_0$. Then
$\gvct{m}=1/|\gvct{q}|^2 \bar{\gvct{q}}\gvct{r}$ is a half integral quaternion,
and since $\gvct{q}$ and $\gvct{r}$ commute it follows from lemma~\ref{sym3}
that $\gvct{q}$ corresponds to a symmetry operation of $\gvct{L}(R(\gvct{r}))$.
Thus we have proved
\begin{lemma}
Let $\gvct{r}=(r_0,r_1,r_2,r_3)$, $\gvct{q}=(1/c) (0,r_1,r_2,r_3)$,
$c=\gcd(r_1,r_2,r_3)$. Then $R(\gvct{q})$ is a symmetry operation of
$\gvct{L}(R(\gvct{r}))$ if and only if $|\gvct{q}|^2$ divides $2r_0$.
\end{lemma}
For the special case $\gvct{r}=(m,n,n,n)$ this lemma states that $(1,1,1)$
is a sixfold axis if and only if $\Sigma(\gvct{r})$ is divisible by $3$,
a result that we have obtained previously by a different method.

\section{Bravais lattices}

So far we have only considered the symmetry of the CSLs, but we can go even
further and compute the Bravais class of the CSL. For some CSLs the Bravais
class follows immediately from theorem~\ref{symcsl}, e.g. for CSLs with
hexagonal and rhombohedral symmetry we know at once that they must belong
to the (unique) hexagonal and rhombohedral Bravais class, respectively.
Consider now an orthorhombic CSL, then we cannot infer from symmetry to which
of the four orthorhombic Bravais classes the CSL belongs. Nevertheless
we can actually compute them.

As an example we consider an orthorhombic CSL generated by
$\gvct{r}=(0,\vec r)$ such that $\vec r=\vec q\times \vec m$, where
$\vec q\cdot\vec m=0$ (case~\ref{orth2} of theorem~\ref{symcsl}).
We assume further
that $\vec r^2$, $\vec q^2$ and $\vec m^2$ are all odd. Then these three
vectors are a basis for the CSL, since $\vec q\cdot \vec r^{(i)}/\vec q^2=-m_i$
and $\vec m\cdot \vec r^{(i)}/\vec m^2=q_i$ are integers and
$|\vec r|\cdot|\vec q|\cdot|\vec m|=\vec r^2=\Sigma(\gvct{r})$. Hence the CSL
is primitive orthorhombic.

Assume now that $\vec r^2$ and $\vec q^2$ are even, and hence $\vec m^2$ odd.
Then $\vec r$, $\vec q$, $\vec m$ span a primitive orthorhombic sublattice of
the CSL of index $2$. Checking all combinations
$\alpha_1\vec r+\alpha_2\vec q+\alpha_3\vec m$, $\alpha_i\in\{0,\pm 1/2\}$,
we find that only the combinations $\pm 1/2\,(\vec r\pm \vec q)$ are integral
vectors, and hence the CSL must be C--face centered orthorhombic, with 
$1/2\,(\vec r\pm \vec q)$, $\vec m$ as a possible basis.

Similarly one can discuss all the other cases. We finally find (for the
conventions of the lattice parameters see~\cite{IT}):
\begin{theo}\label{brav}
Let $\gvct{r}$ be equivalent to $(0,\vec r)$. Then $\gvct{L}(R(\gvct{r}))$
belongs to one of the following Bravais classes (always $0\ne n\ne m\ne 0$):
\begin{enumerate}
\item If $\gvct{r}\sim (m,n,n,n)$ and $3$ does not
divide $\Sigma$, then the CSL is rhombohedral with lattice parameters
$a=\sqrt{2\Sigma}, c=\sqrt{3}$ (triple hexagonal setting).
\item\label{bhex}
If $\gvct{r}\sim (m,n,n,n)$ and $3$ divides $\Sigma$,
then the CSL is hexagonal with lattice parameters
$a=\sqrt{2\Sigma/3}, c=\sqrt{3}$.
\item\label{btetr} If $\gvct{r}\sim (m,n,0,0)$, the CSL is a primitive tetragonal
lattice with lattice parameters $a=\sqrt{\Sigma}, c=1$.
\item\label{borth1}
If $\gvct{r}\sim (m,n,n,0)$,
the CSL is B--face centered orthorhombic with lattice parameters $a=\sqrt{2},
b=\sqrt{\Sigma}, c=\sqrt{2\Sigma}$.
\privateremark{If $m$ is odd $(n,-n,m)$ is parallel to the b-axis and
$(m,-m,-2n)$ is parallel to the c-axis, if $m$ is even and hence $n$ odd,
then $(n,-n,m)$ is parallel to the c-axis and $(m/2,-m/2,-n)$ is parallel to
the b-axis.}
\item\label{borth2}
If $\gvct{r}=(0,\vec r)$ is not equivalent to one of the cases above,
and if there exist two orthogonal integer vectors
$\vec q$ and $\vec m$, $\vec q^2$ and $\vec m^2$ relatively prime, such that
$\vec r=\vec q\times\vec m$ or $\vec r=1/2\ \vec q\times\vec m$, then
the CSL is orthorhombic. In particular:
\begin{enumerate}
\item If $\vec r^2$, $\vec q^2$ and $\vec m^2$ are all odd, the CSL is
primitive orthorhombic with lattice parameters $a=|\vec q|$,
$b=\sqrt{\Sigma}/|\vec q|$, $c=\sqrt{\Sigma}$.
\item If $\vec r^2$ and $\vec q^2$ are even, then the CSL is B--face centered
orthorhombic with lattice parameters $a=|\vec q|$,
$b=\sqrt{2\Sigma}/|\vec q|$, $c=\sqrt{2\Sigma}$.
\privateremark{If $\vec q^2$
is odd and $\vec m^2$ is even, change the roles of $\vec q$ and $\vec m$.}
\item If $\vec r^2$ is odd, $\vec q^2$ and $\vec m^2$ are even
\privateremark{in this case and only in this case
$\vec r=1/2\ \vec q\times\vec m$} the CSL is C--face centered
orthorhombic with lattice parameters $a=|\vec q|$,
$b=2\sqrt{\Sigma}/|\vec q|$, $c=\sqrt{\Sigma}$.
\end{enumerate}
\item\label{bmon}
If $\gvct{r}=(0,\vec r)$ is not equivalent to one of the cases above,
i.e. no decomposition $\vec r=1/2^\ell\ \vec q\times\vec m$ exists, then the
CSL is monoclinic. If $\vec r^2$ is odd, the CSL is primitive monoclinic;
if $\vec r^2$ is even, the CSL is a C-type monoclinic lattice. If $\vec r$
is parallel to the $c$-axis then $c=\sqrt{\Sigma}$
in the first case and $c=\sqrt{2\Sigma}$ in the latter.
\end{enumerate}
\end{theo}

\section{Remarks on non--primitive cubic lattices}

So far, we have only dealt with primitive cubic lattices, but many results
remain true for face--centered and body--centered lattices, too.
The index $\Sigma(R)$ and the number of equivalent and inequivalent CSLs
are the same for all cubic lattices, and so is the symmetry of the
CSLs~\cite{gribo74,baa97,grim76}.

The primitive cubic lattice $\gvct{L}_p$ can be written as a union
$\gvct{L}_p=\gvct{L}_0\cup\gvct{L}_1\cup\gvct{L}_2\cup\gvct{L}_3$ of four
disjoint subsets of $\gvct{L}_p$. Here $\gvct{L}_i$ contains those vectors~$\vec v$
of $\gvct{L}_p$ for which $\vec v^2=i\pmod 4$ holds. Note that only
$\gvct{L}_0=2\gvct{L}_p$ is a lattice. Similarly we can write the body centered
cubic lattice $\gvct{L}_b$ and the face centered cubic lattice $\gvct{L}_f$ as 
$\gvct{L}_b=\gvct{L}_0\cup\gvct{L}_1\cup\gvct{L}_2\cup\gvct{L}_3\cup\frac{1}{2}\gvct{L}_3=
\gvct{L}_p\cup\frac{1}{2}\gvct{L}_3$
and
$\gvct{L}_f=\gvct{L}_0\cup\gvct{L}_1\cup\gvct{L}_2\cup\gvct{L}_3\cup\frac{1}{2}\gvct{L}_2=
\gvct{L}_p\cup\frac{1}{2}\gvct{L}_2$, respectively. One can show that
$\gvct{L}_p(R)\cap\gvct{L}_i=\gvct{L}_i\cap R\gvct{L}_i=:\gvct{L}_i(R)$, and hence
$\gvct{L}_b(R)=\gvct{L}_p(R)\cup\frac{1}{2}\gvct{L}_3(R)$ and
$\gvct{L}_f(R)=\gvct{L}_p(R)\cup\frac{1}{2}\gvct{L}_2(R)$, see e.g.~\cite{gribo74}.
These relations can then be used to prove that $\Sigma(R)$ is the same for
all three types of cubic lattices, that the number of equivalent and
inequivalent CSLs and their symmetry is the same for all cubic lattices as
well.

Of course the lattices itself are different, and so the CSLs usually belong
to different Bravais classes. First we generalize lemma~\ref{cslbasis}
for body and face centered cubic lattices.
\begin{lemma}\label{cslbasisb}
Let $\gvct{L}_b$ be a body centered cubic lattice and
$\gvct{r}=(r_0,r_1,r_2,r_3)$ a primitive quaternion. Then
the CSL $\gvct{L}_b(R(\gvct{r}))$ is the $\Z$--span of the following vectors:
\begin{itemize}
\item  $\vec r^{(0)},\vec r^{(1)},\vec r^{(2)},\vec r^{(3)},
1/2\,(\vec r^{(0)}+\vec r^{(1)}+\vec r^{(2)}+\vec r^{(3)})$
if $|\gvct{r}|^2$ is odd,
\item  $\vec r^{(0)}, 1/2\,(\vec r^{(0)}+\vec r^{(1)}),
1/2\,(\vec r^{(0)}+\vec r^{(2)}), 1/2\,(\vec r^{(0)}+\vec r^{(3)})$ 
if $2\big| |\gvct{r}|^2$ and $4\!\!\not\big|\, |\gvct{r}|^2$,
\item $ 1/2 \vec r^{(0)}, 1/2 \vec r^{(1)}, 1/2 \vec r^{(2)}, 1/2 \vec r^{(3)}$
if $4\big| |\gvct{r}|^2$.
\end{itemize}
\end{lemma}
The situation is a bit nastier for the face centered lattice:
\begin{lemma}\label{cslbasisf}
Let $\gvct{L}_f$ be a face centered cubic lattice and
$\gvct{r}=(r_0,r_1,r_2,r_3)$ a primitive quaternion. Then we have the following
cases:
\begin{itemize}
\item If $\gvct{r}^2=3\pmod 4$ define $\ell_i=1$ if $r_i$ is odd and
$\ell_i=0$ if $r_i$ is even. Then $\gvct{L}_f(R(\gvct{r}))$ is the $\Z$--span of
the vectors
$2^{-\ell_0}\vec r^{(0)}, 2^{-\ell_1}\vec r^{(1)}, 2^{-\ell_2}\vec r^{(2)},
2^{-\ell_3}\vec r^{(3)}$.
\item If $\gvct{r}^2=1\pmod 4$ then $\gvct{L}_f(R(\gvct{r}))$ is the $\Z$--span of
the vectors $\vec r^{(0)},\vec r^{(1)},\vec r^{(2)},\vec r^{(3)}$ and those
combinations $1/2\,(\vec r^{(i)}+\vec r^{(j)})$, for which $r_i+r_j$ is even.
\item If $\gvct{r}^2=2\pmod 4$ define $\ell_i=1$ if $r_i$ is even and
$\ell_i=0$ if $r_i$ is odd. Further define $m_i=0$ if
$1/2\,(r_0-r_1-r_2-r_3)-r_i$ is even and $m_i=1$ if $1/2\,(r_0-r_1-r_2-r_3)-r_i$
is odd. Then $\gvct{L}_f(R(\gvct{r}))$ is the $\Z$--span of
the vectors $2^{-\ell_i}\vec r^{(i)}$ and
$2^{m_i}/4\,(\vec r^{(0)}+\vec r^{(1)}+\vec r^{(2)}+\vec r^{(3)}+2\vec r^{(i)})$,
$i=0,\ldots 3$.
\item If $\gvct{r}^2=0\pmod 4$ then $\gvct{L}_f(R(\gvct{r}))$ is the $\Z$--span of
the vectors $\vec r^{(0)}$,
$1/4\,(\vec r^{(0)}+(-1)^{(r_0-r_i)/2}\vec r^{(i)})$,
$i=1,2,3$.
\end{itemize}
\end{lemma}

In principle one could use these representations to derive the symmetries
and the Bravais class for the CSLs of the body and face centered cubic
lattices. However, it is simpler to derive them from the primitive cubic case
by means of $\gvct{L}_b(R)=\gvct{L}_p(R)\cup\frac{1}{2}\gvct{L}_3(R)$ and
$\gvct{L}_f(R)=\gvct{L}_p(R)\cup\frac{1}{2}\gvct{L}_2(R)$.
Finally we obtain the following results for the Bravais lattices. For
the body centered case we have
\begin{theo}
Let $\gvct{r}$ be equivalent to $(0,\vec r)$. Then $\gvct{L}_b(R(\gvct{r}))$
belongs to one of the following Bravais classes (always $0\ne n\ne m\ne 0$):
\begin{enumerate}
\item If $\gvct{r}\sim (m,n,n,n)$ and $3$ does not
divide $\Sigma$, then the CSL is rhombohedral with lattice parameters
$a=\sqrt{2\Sigma}, c=\sqrt{3}/2$ (triple hexagonal setting).
\item\label{bhexb}
If $\gvct{r}\sim (m,n,n,n)$ and $3$ divides $\Sigma$,
then the CSL is hexagonal with lattice parameters
$a=\sqrt{2\Sigma/3}, c=\sqrt{3}/2$.
\item\label{btetrb} If $\gvct{r}\sim (m,n,0,0)$, the CSL is a body centered
tetragonal lattice
with lattice parameters $a=\sqrt{\Sigma}, c=1$.
\item\label{borth1b}
If $\gvct{r}\sim (m,n,n,0)$ with $m$ odd and $n$ even, then
the CSL is face centered orthorhombic with lattice parameters $a=\sqrt{2},
b=\sqrt{\Sigma}, c=\sqrt{2\Sigma}$. If both $m$ and $n$ are odd then the CSL
is a B--face centered orthorhombic lattice with lattice parameters $a=\sqrt{2},
b=\sqrt{\Sigma}/2, c=\sqrt{2\Sigma}$.
If $m$ is even and not divisible by $4$,
then the CSL is B--face centered orthorhombic with lattice parameters
$a=\sqrt{2}, b=\sqrt{\Sigma}/2, c=\sqrt{2\Sigma}$. If $m$ is divisible by $4$,
then the CSL is face centered orthorhombic with lattice parameters $a=\sqrt{2},
b=\sqrt{\Sigma}, c=\sqrt{2\Sigma}$. 
\privateremark{If both $m$ and $n$ are odd the $b$-axis is parallel to 
$(n,-n,m)$.
If $m$ is even and not divisible by $4$, then the $b$-axis
is parallel to $(m/4,-m/4,-n/2)$.}
\item\label{borth2b}
If $\gvct{r}=(0,\vec r)$ is not equivalent to one of the cases above,
and if there exist two orthogonal integer vectors
$\vec q$ and $\vec m$, $\vec q^2$ and $\vec m^2$ relatively prime, such that
$\vec r=\vec q\times\vec m$ or $\vec r=1/2\ \vec q\times\vec m$, then
the CSL is orthorhombic. In particular:
\begin{enumerate}
\item If $\vec r^2$, $\vec q^2$ and $\vec m^2$ are all odd, the CSL is
body centered orthorhombic with lattice parameters $a=|\vec q|$,
$b=\sqrt{\Sigma}/|\vec q|$, $c=\sqrt{\Sigma}$.
\item If $\vec r^2$ and $\vec q^2$ are even, then the CSL is B--face centered
orthorhombic if $\vec m\in \gvct{L}_3$ and face centered orthorhombic otherwise.
The lattice parameters $a=|\vec q|$,
$b=\sqrt{\Sigma}/2|\vec q|$, $c=\sqrt{2\Sigma}$ and $a=|\vec q|$,
$b=\sqrt{2\Sigma}/|\vec q|$, $c=\sqrt{2\Sigma}$, respectively.
\privateremark{If $\vec q^2$
is odd and $\vec m^2$ is even, change the roles of $\vec q$ and $\vec m$.}
\item If $\vec r^2$ is odd, $\vec q^2$ and $\vec m^2$ are even,
\privateremark{in this case and only in this case
$\vec r=1/2\ \vec q\times\vec m$} the CSL is C--face centered
orthorhombic if $\vec r\in \gvct{L}_3$ and face centered orthorhombic
otherwise. The lattice parameters for these two cases read $a=|\vec q|$,
$b=2\sqrt{\Sigma}/|\vec q|$, $c=\sqrt{\Sigma}/2$ and 
$a=|\vec q|$, $b=2\sqrt{\Sigma}/|\vec q|$,
$c=\sqrt{\Sigma}$,
respectively.
\end{enumerate}
\item\label{bmonb}
If $\gvct{r}=(0,\vec r)$ is not equivalent to one of the cases above,
i.e. no decomposition $\vec r=1/2^\ell\ \vec q\times\vec m$ 
with $\vec q\cdot\vec m=0$ exists then the
CSL is a C-type monoclinic lattice if $\vec r^2\neq 3\pmod 4$. If
$\vec r^2=3\pmod 4$ the CSL is primitive monoclinic.
\end{enumerate}
\end{theo}

Similarly we can determine the Bravais classes in the face centered cubic
case.
\begin{theo}
Let $\gvct{r}$ be equivalent to $(0,\vec r)$. Then $\gvct{L}_f(R(\gvct{r}))$
belongs to one of the following Bravais classes (always $0\ne n\ne m\ne 0$):
\begin{enumerate}
\item If $\gvct{r}\sim (m,n,n,n)$ and $3$ does not
divide $\Sigma$, then the CSL is rhombohedral with lattice parameters
$a=\sqrt{\Sigma/2}, c=\sqrt{3}$ (triple hexagonal setting).
\item\label{bhexf}
If $\gvct{r}\sim (m,n,n,n)$ and $3$ divides $\Sigma$,
then the CSL is hexagonal with lattice parameters
$a=\sqrt{\Sigma/6}, c=\sqrt{3}$.
\item\label{btetrf} If $\gvct{r}\sim (m,n,0,0)$, the CSL is a body centered
tetragonal lattice with lattice parameters $a=\sqrt{\Sigma/2}, c=1$.
\item\label{borth1f}
If $\gvct{r}\sim (m,n,n,0)$ with $m$ odd and $n$ even, then
the CSL is body centered orthorhombic with lattice parameters $a=1/\sqrt{2},
b=\sqrt{\Sigma/2}, c=\sqrt{\Sigma}$. If both $m$ and $n$ are odd, then
the CSL is a C--face centered orthorhombic lattice
with lattice parameters $a=1/\sqrt{2}, b=\sqrt{\Sigma/2}, c=\sqrt{\Sigma}$. 
If $m$ is even and not divisible by $4$ then the CSL is C--face centered
orthorhombic lattice
with lattice parameters $a=1/\sqrt{2},
b=\sqrt{\Sigma/2}, c=\sqrt{\Sigma}$. If $m$ is divisible by $4$ the CSL
is again a body centered orthorhombic lattice
with lattice parameters $a=1/\sqrt{2},
b=\sqrt{\Sigma/2}, c=\sqrt{\Sigma}$.
\privateremark{If both $m$ and $n$ are odd the $b$-axis is parallel to
$(m/2,-m/2,-n)$. The same is true for $m$ even and not divisible by $4$.}
\item\label{borth2f}
If $\gvct{r}=(0,\vec r)$ is not equivalent to one of the cases above,
and if there exist two orthogonal integer vectors
$\vec q$ and $\vec m$, $\vec q^2$ and $\vec m^2$ relatively prime, such that
$\vec r=\vec q\times\vec m$ or $\vec r=1/2\ \vec q\times\vec m$, then
the CSL is orthorhombic. In particular:
\begin{enumerate}
\item If $\vec r^2$, $\vec q^2$ and $\vec m^2$ are all odd, the CSL is
face centered orthorhombic with lattice parameters $a=|\vec q|$,
$b=\sqrt{\Sigma}/|\vec q|$, $c=\sqrt{\Sigma}$. 
\item If $\vec r^2$ and $\vec q^2$ are even, then the CSL is B--face centered
orthorhombic with lattice parameters $a=|\vec q|/2$,
$b=\sqrt{2\Sigma}/|\vec q|$, $c=\sqrt{\Sigma/2}$ if
$1/2\,(\vec r+\vec q)\in\gvct{L}_2$. If $1/2\,(\vec r+\vec q)\not\in\gvct{L}_2$
then $1/2\,(\vec r+\vec q)+\vec m\in\gvct{L}_2$ and the CSL is body centered
orthorhombic with lattice parameters $a=|\vec q|/2$,
$b=\sqrt{2\Sigma}/|\vec q|$, $c=\sqrt{\Sigma/2}$.
\privateremark{If $\vec q^2$
is odd and $\vec m^2$ is even, change the roles of $\vec q$ and $\vec m$.}
\item If $\vec r^2$ is odd, $\vec q^2$ and $\vec m^2$ are even
\privateremark{in this case and only in this case
$\vec r=1/2\ \vec q\times\vec m$} then the CSL is C--face centered
orthorhombic or body centered orthorhombic, according to whether
$1/2\,(\vec q+\vec m)\in\gvct{L}_2$ or $1/2\,(\vec q+\vec m)+\vec r\in\gvct{L}_2$.
In both cases the lattice parameters are $a=|\vec q|/2$,
$b=\sqrt{\Sigma}/|\vec q|$, $c=\sqrt{\Sigma}$.
\end{enumerate}
\item\label{bmonf}
If $\gvct{r}=(0,\vec r)$ is not equivalent to one of the cases above,
i.e. no decomposition $\vec r=1/2^\ell\ \vec q\times\vec m$
with $\vec q\cdot\vec m=0$ exists then the
CSL is monoclinic. It is centered monoclinic except if $\vec r^2=3\pmod 4$,
where it is primitive monoclinic.
\end{enumerate}
\end{theo}

Using these theorems one can immediately determine the Bravais class for each
CSL $\gvct{L}(R)$. Table~\ref{tabbrav} lists them for all CSLs
with $\Sigma\leq 59$. 

\begin{table}
\begin{tabular}{|r|l|r|r|r|}
\hline
$\Sigma$ & $\gvct{r}$ & CSL ($cP$) & CSL ($cI$) & CSL ($cF$)\\
\hline
$3$ & $(0,1,1,1)$ & $hP$ & $hP$ & $hP$ \\
$5$ & $(2,1,0,0)$ & $tP$ & $tI$ & $tI$ \\
$7$ & $(2,1,1,1)$ & $hR$ & $hR$ & $hR$\\
$9$ & $(1,2,2,0)$ & $oC$ & $oF$ & $oI$\\
$11$ & $(3,1,1,0)$ & $oC$ & $oC$ & $oC$\\
$13$ & $(1,2,2,2)$ & $hR$ & $hR$ & $hR$\\
$13$ & $(3,2,0,0)$ & $tP$ & $tI$ & $tI$ \\
$15$ & $(0,5,2,1)$
 & $oC$ & $oF$ & $oI$\\
$17$ & $(4,1,0,0)$ & $tP$ & $tI$ & $tI$ \\
$17$ & $(3,2,2,0)$ & $oC$ & $oF$ & $oI$\\
$19$ & $(4,1,1,1)$ & $hR$ & $hR$ & $hR$\\
$19$ & $(1,3,3,0)$ & $oC$ & $oC$ & $oC$\\
$21$ & $(3,2,2,2)$ & $hP$ & $hP$ & $hP$ \\
$21$ & $(0,4,2,1)$
 & $oC$ &  $oF$ & $oI$\\
$23$ & $(0,6,3,1)$ & $mC$ & $mC$ & $mC$\\
$25$ & $(4,3,0,0)$ & $tP$ & $tI$ & $tI$ \\
$25$ & $(0,5,4,3)$ & $mC$ & $mC$ & $mC$\\
$27$ & $(5,1,1,0)$ & $oC$ & $oC$ & $oC$\\
$27$ & $(0,7,2,1)$ & $mC$ & $mC$ & $mC$\\
$29$ & $(5,2,0,0)$ & $tP$ & $tI$ & $tI$ \\
$29$ & $(0,4,3,2)$ & $mP$ & $mC$ & $mC$\\
$31$ & $(2,3,3,3)$ & $hR$ & $hR$ & $hR$\\
$31$ & $(0,7,3,2)$ & $mC$ & $mC$ & $mC$\\
$33$ & $(5,2,2,0)$ & $oC$ & $oF$ & $oI$\\
$33$ & $(1,4,4,0)$ & $oC$ & $oF$ & $oI$\\
$33$ & $(0,7,4,1)$
 & $oC$ & $oC$ & $oC$\\
$35$ & $(0,5,3,1)$
 & $oC$ & $oC$ & $oC$\\
$35$ & $(0,6,5,3)$
 & $oC$ & $oF$ & $oI$\\
$37$ & $(6,1,0,0)$ & $tP$ & $tI$ & $tI$ \\
$37$ & $(5,2,2,2)$ & $hR$ & $hR$ & $hR$\\
$37$ & $(0,8,3,1)$ & $mC$ & $mC$ & $mC$\\
$39$ & $(6,1,1,1)$ & $hP$ & $hP$ & $hP$ \\
$39$ & \begin{minipage}{1.75cm}\noindent$(5,3,2,1)$ $(-5,3,2,1)$\end{minipage}
 & $mC$ & $mC$ & $mC$\\
$41$ & $(5,4,0,0)$ & $tP$ & $tI$ & $tI$ \\
$41$ & $(3,4,4,0)$ & $oC$ & $oF$ & $oI$\\
$41$ & $(0,6,2,1)$ & $mP$ & $mC$ & $mC$\\
$43$ & $(4,3,3,3)$ & $hR$ & $hR$ & $hR$\\
$43$ & $(5,3,3,0)$ & $oC$ & $oC$ & $oC$\\
$43$ & $(0,9,2,1)$ & $mC$ & $mC$ & $mC$\\
$45$ & $(0,5,4,2)$
 & $oP$ & $oI$ & $oF$\\
$45$ & $(0,8,5,1)$ & $oC$ & $oC$ & $oC$\\
$45$ & $(0,7,5,4)$ & $oC$ & $oC$ & $oC$\\
$47$ & $(0,9,3,2)$ & $mC$ & $mC$ & $mC$\\
$47$ & $(0,7,6,3)$ & $mC$ & $mC$ & $mC$\\
$49$ & $(1,4,4,4)$ & $hR$ & $hR$ & $hR$\\
$49$ & $(0,6,3,2)$ & $mP$ & $mC$ & $mC$\\
$49$ & $(0,9,4,1)$ & $mC$ & $mC$ & $mC$\\
$51$ & $(7,1,1,0)$ & $oC$ & $oC$ & $oC$\\
$51$ & $(1,5,5,0)$ & $oC$ & $oC$ & $oC$\\
$51$ & \begin{minipage}{1.75cm}\noindent $(5,4,3,1)$ $(-5,4,3,1)$\end{minipage}
 & $mP$ & $mC$ & $mC$\\
$53$ & $(7,2,0,0)$ & $tP$ & $tI$ & $tI$ \\
$53$ & $(0,6,4,1)$ & $mP$ & $mC$ & $mC$\\
$53$ & $(0,9,4,3)$ & $mC$ & $mC$ & $mC$\\
$55$ & $(0,10,3,1)$
 & $oC$ & $oC$ & $oC$\\
$55$ & $(0,9,5,2)$ & $mC$ & $mC$ & $mC$\\
$55$ & $(0,7,6,5)$ & $mC$ & $mC$ & $mC$\\
\hline
\end{tabular}
\end{table}
\begin{table}
\begin{tabular}{|r|l|r|r|r|}
\hline
$\Sigma$ & $\gvct{r}$ & CSL (cP) & CSL (cI) & CSL (cF)\\
\hline
$57$ & $(3,4,4,4)$ & $hP$ & $hP$ & $hP$ \\
$57$ & $(7,2,2,0)$ & $oC$ & $oF$ & $oI$\\
$57$ & $(5,4,4,0)$ & $oC$ & $oF$ & $oI$\\
$57$ & \begin{minipage}{1.75cm}\noindent $(6,4,2,1)$ $(-6,4,2,1)$\end{minipage}
 & $mC$ & $mC$ & $mC$\\
$59$ & $(3,5,5,0)$ & $oC$ & $oC$ & $oC$\\
$59$ & $(0,7,3,1)$ & $mP$ & $mP$ & $mP$\\
$59$ & $(0,9,6,1)$ & $mC$ & $mC$ & $mC$\\
\hline
\end{tabular}
\caption{\label{tabbrav}Bravais classes of the CSLs with $\Sigma\leq 59$}
\end{table}

\section{Further remarks and outlook}

We have derived the symmetry properties of the CSLs by making intensive
use of quaternions and the proofs are mainly algebraic. A crystallographer
not familiar with quaternions might be interested in a more geometric
development of this topic. Indeed, one can prove most theorems with
geometrical methods. We briefly sketch how this can be done for the primitive
cubic case. Let $R=R(\gvct{r})=R(r_0,\vec r)$ If $\vec v\in\gvct{L}(R)$ then
\BAL
R^{-1}\vec v=\frac{1}{r_0^2+\vec r^2}\left((r_0^2-\vec r^2)\vec v
-2r_0 \vec r\times\vec v + 2(\vec r\cdot\vec v)\vec r \right)\in \gvct{L},
\EAL
i.e. $R^{-1}\vec v$ must be an integer vector. This expression simplifies
if $R$ is a rotation through $\pi$, where we have $r_0=0$. It then
follows that $R^{-1}\vec v\in\gvct{L}$ if and only if $\vec r^2$ divides
$2\vec r\cdot\vec v$. If $Q$ is a symmetry operation of $\gvct{L}(R)$,
then $Q^{-1}\vec v\in\gvct{L}$ and $R^{-1}Q^{-1}\vec v\in\gvct{L}$ for all
$\vec v\in\gvct{L}$. If we assume that $Q$ is a rotation through $\pi$ around
the axis $\vec q$ we get the following two conditions: $\vec q^2$ must divide
$2\vec q\cdot\vec v$ and
\BAL
\frac{2\vec q\cdot\vec v}{\vec q^2}\left(-\vec q+
\frac{2\vec r\cdot\vec v}{\vec r^2}\vec r \right)\in \gvct{L},
\EAL
which implies that
$\frac{2\vec q\cdot\vec v}{\vec q^2}\frac{2\vec r\cdot\vec v}{\vec r^2}$ is
an integer. One shows further that there exist an integer $n$ and an integer
vector $\vec c$ orthogonal to $\vec r$
such that $\vec r=n\vec q+\vec c$ if $\vec q^2$ is odd and
$\vec r=n/2 \vec q + 1/2 \vec c$ if $\vec q^2$ is even. Now one can
prove that $\vec r^2$ must divide $4\vec c^2$.
If $\vec q^2$ is even we have the
stricter condition that $\vec r^2$ must divide $\vec c^2$. These conditions
limit the possible values of $n$ and $\vec c$. If one checks all the possible
cases (which is a bit tedious) one finally arrives at theorem~\ref{brav}.

We have answered the question which symmetries a CSL has and to which
Bravais class it belongs for all $\gvct{r}$ equivalent to $(0,\vec r)$,
It would be interesting to answer the question also for the general case.

We have shown that a CSL has orthorhombic symmetry if $\gvct{r}$ can be written
as $2^\ell\gvct{r}=\gvct{q}\gvct{m}=-\gvct{m}\gvct{q}$, or equivalently
$\vec r=\vec q\times \vec m$, where $\vec q\cdot\vec m=0$. Here the question
arises under what conditions such a decomposition exists, and how many
inequivalent decompositions exist for a given $\Sigma=2^{-\ell}|\gvct{r}|^2$.
For a fixed $\gvct{r}$, such a decomposition is unique up to sign changes
and permutations, unless $\gvct{r}\sim(0,1,1,1)$ or $(0,1,0,0)$. 
This question is related
to the number of inequivalent but congruent CSLs. Grimmer suggests
a formula relating the number of congruent CSLs with several other properties
like the symmetry of the CSL\cite{grim76}. This formula is based on the
analysis of the CSLs up to $\Sigma=199$. Using our results, this formula can
be proved for the hexagonal, tetragonal and rhombohedral CSLs. For the
orthorhombic case, one would need a formula for the number of representations
of the kind $2^m\gvct{r}=\gvct{q}\gvct{m}=-\gvct{m}\gvct{q}$.

\section*{Acknowledgements}

The author is very grateful to Michael Baake for interesting discussions
on the present subject. Financial support by the Austrian
Academy of Sciences (APART-program) is gratefully acknowledged.


\end{document}